\documentclass[11pt,draft]{article}
\usepackage[cp1250]{inputenc}

\usepackage{amssymb}
\usepackage{amsmath}
\usepackage{fancyhdr}
\usepackage{amsthm}

\title{On identities in centrally nilpotent Moufang loops
and centrally nilpotent $A$-loops}
\date{}
\author{N. I. Sandu}

\begin{document}

 \maketitle

\begin{abstract}
This paper  proves that the variety generated by a centrally
nilpotent Moufang loop (or centrally nilpotent $A$-loop) is
finitely based.
\smallskip\\
\textbf{Key words}: Moufang loop, $A$-loop, centrally nilpotent
loop, lower central series, commutator-associator  of type
$(\alpha, \beta)$ of weight $n$  com\-mutator-associators of type
$(\mu)$ of weight $n$), fully invariant  subloop, word subloop,
basis of identities.

\textbf{Mathematics of subject classification}: 20N05.
\end{abstract}

Combining the  techniques used by Lyndon \cite{Lyn} in showing
that the identities of a nilpotent group are finitely based with
some related ideas of Higman \cite{Hig} and the results of Bruck
\cite{Br58} on the structure of commutative Moufang loops Evans
\cite{Evans} prove that the identities of centrally nilpotent
commutative Moufang loop are finitely based. In this paper the
Evans's result  is reinforced for centrally nilpotent Moufang
loops and centrally nilpotent $A$-loops (Theorem 5.4), following
the ideas of \cite{Evans}.

In  \cite{Br46}, \cite{Br58} Bruck developed and investigated in
details the theory of nilpotency, particularly centrally
nilpotency, for loops on the basis of  nilpotency for groups.

The theory of centrally nilpotent loops essentially differs from
the theory of nilpotent groups. In some of his papers, Sandu
offers methods to research   centrally nilpotent loops, distinct
from \cite{Br46}, \cite{Br58}. They are  systematized in
\cite{CovSan1} and \cite{CovSan2}. The  notions of commutators, of
commutators of weight $n$, of lower central series, of upper
central series etc. for groups are replaced by the notions of
commutator-associator, of commutator-associators of weight $n$, of
lower central series, of upper central series, etc. for loops
respectively.   The basic result of this paper, Theorem 5.4, is
proved using the mentioned replacement. Particularly, it  uses the
structural construction of  members of lower central series with
the help of commutator-associators of type $(\alpha, \beta)$.

Theorem 5.4 from this paper  is related to the main results of
papers \cite{Ur2000}, \cite{Ur2010}. But, unfortunately, these
results from \cite{Ur2000}, \cite{Ur2010}  cannot be  consider as
proved. They are described in details in \cite{Sandu14}.

We will use the notation  $\Box$ to mark the completion of a
proof.

\section{Preliminaries}

Let us remind some notions and results from the loop theory, which
can be found in \cite{Br58} (see, also \cite{Bel},
\cite{Smith}).\vspace*{0.1cm}

A \textit{quasigroup} is a non-empty set $Q$ together with a
binary operation $Q \times Q \rightarrow Q; (x, y) \rightarrow xy$
such that the equations $$ax = c, \quad yb = c \eqno{(1.1)}$$ have
unique solutions. A \textit{loop}  $(Q, \cdot, 1)$ is a quasigroup
$(Q, \cdot)$  with a base point, or distinguished element, $1 \in
Q$  satisfying the equations $1x = x1 = x$ for all $x \in
Q$.\vspace*{0.1cm}

The \textit{multiplication group} $ \frak C(Q)$ of the arbitrary
loop $Q = (Q,\cdot, 1)$ is generated by all mappings $R(x)$,
$L(x)$, where $R(x)y = yx, L(x)y = xy$, of the loop $Q$. The
\textit{inner mapping group} $\frak I(Q)$ of $Q$ is the subgroup
of $\frak C(Q)$, generated by the mappings $T(x), R(x,y), L(x,y)$,
for all $x, y$ in $Q$, where $$ T(x) = L^{-1}(x)R(x), \quad R(x,y)
= R^{-1}(xy)R(y)R(x),$$ $$L(x,y) = L^{-1}(xy)L(x)L(y).
\eqno{(1.2)}. $$

For the loop $Q$ with identity $1$ $$\frak I(Q) = \{\alpha \in
\frak C(Q) \vert \alpha 1 = 1 \}. \eqno{(1.3)}$$

The subloop $H$ of the loop $Q$ is called \textit{normal}
(\textit{invariant}) in $Q$, if
$$ xH = Hx, \quad x\cdot yH = xy\cdot H, \quad H\cdot xy = Hx\cdot
y \eqno{(1.4)} $$ or by (1.2) $$ T(x)H = H, \quad L(x,y)H = H,
\quad R(x,y)H = H \eqno{(1.5)} $$ for every $x, y \in
Q$.\vspace*{0.1cm}

We will use the notation $<M>$ for the subloop of the loop $Q$
generated by set $M \subseteq Q$. Let $H, K$ be a subloop of the
loop $Q$ such that $K$ is normal in $<H, K>$. Then $<H, K> = HK =
KH$. \vspace*{0.1cm}

\textbf{Lemma 1.1.} \textit{Let $H$ be a normal subloop of the
loop $Q$. If $H$ is generated as a normal subloop by set $S
\subseteq Q$, then $H$, as a subloop, is generated by set
$\{\varphi s \vert s \in S, \varphi \in \frak I(Q)\}$, where
$\frak I(Q)$ is the inner mapping group of $Q$.} \vspace*{0.1cm}

For arbitrary elements $x, y, z$ of the loop $Q$ the
\textit{commutator}, $(x,y)$, and \textit{associator}, $(x,y,z)$,
are defined by $$xy = (yx)[x,y], \quad xy\cdot z = (x\cdot
yz)[x,y,z]. \eqno{(1.6)}$$

The \textit{left nucleus}, $N_{\lambda}(Q)$, of the loop $Q$ is
the associative subloop $N_{\lambda}(Q) = \{a \in Q \vert [a,x,y]
= 1 \quad \forall x, y \in Q\}$, the \textit{middle nucleus},
$N_{\mu}(Q)$, of $Q$ is the associative subloop $N_{\mu}(Q) = \{a
\in Q \vert [x,a,y] = 1 \quad \forall x, y \in Q\}$ and the
\textit{right nucleus}, $N_{\rho}(Q)$, of $Q$ is the associative
subloop $N_{\rho}(Q) = \{a \in Q \vert [x,y,a] = 1 \quad \forall
x, y \in Q\}$.\vspace*{0.1cm}

The \textit{nucleus}, $N(Q)$, is defined by $N(Q) = N_{\lambda}(Q)
\cap N_{\mu}(Q) \cap N_{\rho}(Q)$ and the \textit{centre}, $Z(Q)$,
of the loop $Q$ is a normal associative and commutative subloop
$Z(Q) = \{a \in N(Q) \vert [a,x] = 1 \quad \forall x \in Q\}$.
Moreover, every subloop $H \subseteq Z(Q)$ is normal in $Q$.
\vspace*{0.1cm}

 The (transfinite) \textit{upper central series},
$\{Z_\alpha\}$, of a loop $Q$  is defined inductively as follows:

i) $Z_0 = 1$;

ii) for any ordinal $\alpha$ $Z_{\alpha + 1}$ is the unique
subloop of $Q$ such that \break $Z_{\alpha + 1}/Z_{\alpha} =
Z(Q/Z_{\alpha})$;

iii) if $\alpha$ is a limit ordinal, $Z_\alpha$ is the union of
all $Z_\beta$ with $\beta < \alpha$.\vspace*{0.1cm}

For any normal subloop of a loop $Q$ we define $(N,Q)$ as the
subloop $A$, whose existence is guaranteed by validity of
condition: if $A$ is the intersection of all normal subloops $K$
of $Q$ such that $NK/K$ is a subloop of $Z(Q/K)$, then $NA/A$ is a
subloop of $Z(Q/A)$.

The (transfinite) \textit{lower central series}, $\{Q_\alpha\}$,
of the loop $Q$  is defined inductively as follows:

i) $Q_0 = Q$;

ii) for any ordinal $\alpha$ $Q_{\alpha + 1} = (Q_\alpha, Q)$;

iii) if $\alpha$ is a limit ordinal, $Q_\alpha$ is the
intersection of all $Q_\beta$ with $\beta <
\alpha$.\vspace*{0.1cm}

Clearly, $Q_{\alpha}$,  $Z_{\alpha}$  are normal subloops of
$Q$.\vspace*{0.1cm}

The loop $Q$ is called \textit{transfinitely centrally nilpotent}
(respect. \textit{transfinitely upper centrally nilpotent}   if
$Q_{\lambda} = 1$ (respect. $Z_{\lambda} = Q$) and
\textit{centrally nilpotent} (respect. \textit{upper centrally
nilpotent} if, in addition, $\lambda$ is finite. The smallest
ordinal $\lambda$ such that $Q_{\lambda} = 1$, is called
\textit{centrally nilpotence class.}\vspace*{0.1cm}

The loop is \textit{Moufang} if it satisfies the equivalent
identities: $$x(y\cdot xz) = (xy\cdot x)z; \quad (xy\cdot z)y =
x(y\cdot zy); \quad xy\cdot zx = (xy\cdot z)x.\eqno{(1.7)}$$

The Moufang loop is diassiciative, i.e. every two elements
generate a subgroup.\vspace*{0.1cm}

The Moufang loop is an $IP$-loop. Then it satisfies the
identities: $$x^{-1}\cdot xy = y; \quad yx\cdot x^{-1} = y; \quad
(xy)^{-1} = y^{-1}x^{-1}; \quad (x^{-1})^{-1} = x.\eqno{(1.8)}$$

\textbf{Lemma 1.2.} \textit{In a Moufang loop $G$, the equation
$[a,b,c] = 1$ implies each of the equations obtained by permuting
$a, b, c$ or replacing any of these elements by their
inverses.}\vspace*{0.1cm}

\textbf{Lemma 1.3.} \textit{Let $a, b, c, d$ be four elements of
the Moufang loop $G$ each three of which associate (satisfy $[x,
y, z] = 1$). Then the following equations are equivalent; (i)
$[ab, c, d] = 1$; (ii) $[cd, a, b] = 1$; $[bc, d, a] =
1$.}\vspace*{0.1cm}

\textbf{Lemma 1.4} \textit{Let $G$ be a Moufang loop. Then $G$
satisfies all or none of the identities:} (i) $[[x,y,z],x] = 1$;
(ii) $[x,y,[y,z]] = 1$; (iii) $[x,y,z]^{-1} = [x^{-1},y,z]$; (iv)
$[x,y,z] = [x,z,y^{-1}]$. \textit{When these identities hold, the
following identities hold:}
$$[x,y,z] = [y,z,x] = [y,x,z]^{-1}, \eqno{(1.9)}$$
$$[xy,z] = [x,z][[x,z],y][y,z][x,y,z]^3,\eqno{(1.10)}$$
\vspace*{0.1cm}

In a Moufang loop the following identities hold:
$$L(z.y)x = x[x,y,z]^{-1}; \eqno{(1.11)}$$
$$L(z,u)(xy)\cdot [z^{-1},u^{-1}] = L(z,u)x \cdot (y\cdot
L(z,u)[z^{-1},u^{-1}]). \eqno{(1.12)}$$

Let $(Q,\cdot,1)$ be a  $IP$-loop. Then  $N(Q) = N_{\lambda}(Q) =
N_{\mu}(Q) = N_{\rho}(Q)$ and  $Z(Q) = \{a \in N(Q) \vert [a,x] =
1 \quad \forall x \in Q\}$.\vspace*{0.1cm}

R. H. Bruck and L. J. Paige \cite{BP} defined an \textit{$A$-loop}
to be a loop in which every inner mapping is an automorphism. Many
of the basic theorems about $A$-loops are contained in \cite{BP};
for example, \textit{$A$-loops are always power associative}
(every $<x>$ is a group), but not necessarily diassociative.

The following statements show an essential difference between
$A$-loops and Moufang loops.\vspace*{0.1cm}

\textbf{Proposition 1.5} \cite{KKP}. \textit{For an $A$-loop $L$,
the following are equivalent:}

\textit{1. $L$ is $IP$-loop;}

\textit{2. $L$ has the alternative property, i.e. $x(xy) = x^2y$
and $(yx)x = yx^2$ for all $x, y \in L$;}

\textit{3. $L$ is diassociative;}

\textit{4. $L$ is a Moufang loop.}

\textbf{Lemma 1.6} \cite{Br58}, \cite{Bel}. \textit{For any
$IP$-loop $Q$ $N_{\lambda}(Q) = N_{\rho}(Q) =
N_{\mu}(Q)$.}\vspace*{0.1cm}

\textbf{Lemma 1.7} \cite{BP}. \textit{For any $A$-loop $Q$
$N_{\lambda}(Q) = N_{\rho}(Q) \subseteq N_{\mu}(Q)$.}

\section{Centrally nilpotent  loops}

Let $Q$ be a loop. The series of normal subloops $Q = C_0
\supseteq C_1 \supseteq \ldots \supseteq C_r = 1$ of loop $Q$ is
called \textit{centrally nilpotent}, if
$$(C_i,Q) \subseteq C_{i+1} \quad \text{for all} \quad i.
\eqno{(2.1)}$$ or, equivalently,
$$C_i/C_{i+1} \subseteq Z(Q/C_{i+1}) \quad\text{for all}\quad i \eqno{(2.2)}$$

\textbf{Lemma 2.1.} \textit{Let $\{C_i\}$ be a centrally nilpotent
series, $\{Z_i\}$ be the upper centrally nilpotent series,
$\{Q_i\}$ be the lower centrally nilpotent series of the loop $Q$.
Then $C_{r - i} \subseteq Z_i$, $C_i \supseteq Q_i$, for $i = 0,
1,\ldots,r$.} \vspace*{0.1cm}

\textbf{Proof.} We have $C_0 = Q = Q_0$. Assume that $C_i
\supseteq Q_i$. By (2.1) $(C_i,Q) \subseteq C_{i+1}$. But then
$Q_{i+1} = (Q_i,Q) \subseteq (C_i,Q) \subseteq C_{i+1}$. We assume
now that $C_{r - i} \subseteq Z_i$ for a certain $i$. Then the
loop $Q/Z_i$ is the homomorphic image of the loop $Q/C_{r - i}$
with kernel $Z_i/C_{r - i}$. But by (2.2)
$$C_{r - i - 1}/C_{r - i} \subseteq Z(Q/C_{r - i}),$$
from where it follows that the homomorphic image of subloop $C_{r
- i - 1}/C_{r - i}$ must lie in the centre $Z(Q/Z_i)$. It is clear
that this image is the subloop $(C_{r - i - 1} \cup Z_i)/Z_i$,
while $Z(Q/Z_i) = Z_{i + 1}/Z_i$. Consequently, $C_{r - i - 1}
\subseteq C_{r-i}\cup Z_i \subseteq Z_{i-1}$.
$\Box$\vspace*{0.1cm}

\textbf{Proposition 2.2.}  \textit{A loop $Q$ is centrally
nilpotent of class $n$ if and only if its upper and lower
centrally nilpotent series have respectively the form}
$$1 = Z_0 \subset Z_1 \subset \ldots \subset Z_n = Q,
\quad Q = Q_0 \supset Q_1 \supset \ldots \supset Q_n = 1.$$

\textbf{Proof.} The statement of the ptoposition for upper
centrally nilpotent series results from the definition of
centrally nilpotent loop. Further, if an centrally nilpotent
series of  length $n$ exists, then from Lemma 2.1 it follows that
the length of the upper and lower central series do not exceed
$n$. But, as there is a term by term inclusion between the
elements of these series, their lengths are  equal, and the series
have the indicated form. $\Box$\vspace*{0.1cm}

We follow \cite{CovSan1}. Let $Q$ be a loop and let $a, b, c \in
Q$. We denote the solution of equation $ab \cdot c = ax \cdot bc$
(respect. $c \cdot ba = cb \cdot xa$) by $\alpha(a,b,c)$ (respect.
$\beta(a, b, c)$) and call it the \textit{associator of type
$\alpha$} (respect. \textit{of type $\beta$}) of elements $a, b,
c$.\vspace*{0.1cm}

The \textit{commutator} of elements $a, b \in Q$ $(a,b)$ is
defined by the equality $ab = b(a(a,b))$.\vspace*{0.1cm}

The last definitions will be written in the form
$$T(b)a = a(a,b), R(b,c)a = a\alpha(a,b,c), L(c,b)a = \beta(a,b,c)a \eqno{(2.3)}$$
if we use (1.2) (see, also, \cite{Br58}).\vspace*{0.1cm}

\textbf{Lemma 2.3} \textit{Any $IP$-loop satisfies the identity}
$$\alpha(x,y,z)^{-1} = \beta(x^{-1},y^{-1},z^{-1}).$$

\textbf{Proof.} From the identity $xy \cdot z =
x\alpha(x,y,z)\cdot yz$ by (1.8) we get $(xy \cdot z)^{-1} =
(x\alpha (x,y,z) \cdot yz)^{-1}$, $z^{-1} \cdot y^{-1}x^{-1} =
z^{-1}y^{-1} \cdot \alpha(x,y,z)^{-1}x^{-1}$ and from the
definition of the associator of type $\beta$ it follows that
$\alpha(x,y,z)^{-1} = \break \beta(x^{-1},y^{-1},z^{-1})$.
$\Box$\vspace*{0.1cm}

\textbf{Lemma 2.4.} \textit{Any Moufang loop $Q$ satisfies the
identities $[x,y,z]^{-1} = \alpha(x,z^{-1},y^{-1})$, $[x,y,z] =
\beta(x^{-1},z,y)$, $[x,y] = (x,y)$.}\vspace*{0.1cm}

\textbf{Proof.} Let $a, b, c \in Q$. By (1.6) $ab\cdot c = (a\cdot
bc)[a,b,c], (ab\cdot c)[a,b,c]^{-1}$
$$= a\cdot bc, [a,b,c]^{-1} = (ab\cdot c)^{-1}(a\cdot bc),$$
$$[a,b,c]^{-1} = (c^{-1}\cdot b^{-1}a^{-1})(a\cdot bc)$$
and by the first identity from (1.7) we have
$$a[a,b,c]^{-1} = a((c^{-1}\cdot b^{-1}a^{-1})(a\cdot bc)) =$$
$$(a(c^{-1}\cdot b^{-1}a^{-1})\cdot a)(bc) = ((ac^{-1})(b^{-1}a^{-1}\cdot a))
(bc) = (ac^{-1}\cdot b^{-1})(bc).$$ Hence $a[a,b,c]^{-1}\cdot
c^{-1}b^{-1} = ac^{-1}\cdot b^{-1}$. Further,
$$(a[a,b,c]^{-1}\cdot c^{-1}b^{-1})^{-1} = (ac^{-1}\cdot b^{-1})^{-1},
bc\cdot [a,b,c]a^{-1} = b\cdot ca^{-1}.$$ Then, from here and from
the definition of the associator of type $\alpha$ or $\beta$, it
follows that $[a,b,c]^{-1} = \alpha(a,c^{-1},b^{-1}), \quad
[a,b,c] = \beta(a^{-1},c,b)$.

Any Moufang loop is diassociative, so $[a,b] = (a,b)$.
$\Box$\vspace*{0.1cm}

Let $Q$ be a loop and let $A, B, C$ be non-empty subsets of $Q$.
We denote $\alpha(A,B,C) = <\alpha(a,b,c) \vert a \in A, b \in B,
c \in C>$, $\beta(A,B,C) = <\beta(a,b,c) \vert a \in A, b \in B, c
\in C>$, $[A,B,C] = <[a,b,c] \vert a \in A, b \in B, c \in C>$,
$(A,B) = <(a,b) \vert a \in A, b \in B>$, $[A,B] = <[a,b] \vert a
\in A, b \in B>$.\vspace*{0.1cm}

\textbf{Lemma 2.5.} \textit{Let $N$ be a normal subloop of the
loop $Q$. Then the subloop $H$, generated by the set
$\alpha(N,Q,Q) \cup \beta(N,Q,Q) \cup(N,Q)$ is normal in
$Q$.}\vspace*{0.1cm}

\textbf{Proof.} Let $n \in N$, $x, y \in Q$. From (2.3) we get
$$\alpha(n,x,y) = L^{-1}(n)R(x,y)n,$$
$$\beta(n,x,y) = R^{-1}(n)L(x,y)n, (n,x) = L^{-1}(n)T(x)n. \eqno{(2.4)}$$
The subloop $N$ is normal in $Q$, then by (1.5) $\alpha(n,x,y),
\beta(n,x,y), (n,x) \in N$ for all $n \in N$ and all $x, y \in Q$.
Then $H \subseteq N$. Let $h \in H$. According to (2.3) we get
$T(x)h = h(h,x) \in H$, $L(x,y)h = \beta(h,y,x) \in H$, $R(x,y)h =
h\alpha(h,x,y) \in H$. Hence by (1.5) the subloop $H$ is normal in
$Q$.$\Box$\vspace*{0.1cm}

\textbf{Proposition 2.6.} \textit{A subloop $H$ of a loop  $Q$ is
normal in $Q$ if and only if (i) $\alpha(H,Q,Q) \subseteq H$,
$\beta(H,Q,Q) \subseteq H$, $(H,Q) \subseteq H$.}\vspace*{0.1cm}

\textbf{Proof.} Let $H$ be a normal subloop of $Q$. Then from
(1.5), (2.3) and  Lemma  2.4 it follows that $H$ satisfies the
inclusions of Proposition 2.6. The inverse statement of
Proposition 2.6 follows from Lemma 2.5. $\Box$\vspace*{0.1cm}

Let $H$ be a normal subloop of a loop $Q$. We denote
$\mathcal{Z}_H(Q) = \{a\in Q \vert \alpha(a,Q,Q) \subseteq H,
\beta(a,Q,Q) \subseteq H, (a, Q) \subseteq H\}$.

Let $E = \{1\}$ and let $\mathcal{Z}_E(Q) = \mathcal{Z}(Q)$. We
prove that
$$\mathcal{Z}(Q) = Z(Q), \eqno{(2.5)}$$
where $Z(Q)$ means the centre of loop $Q$. Really, from the
definitions of the associators $\alpha(x,y,z)$, $\beta(x,y,z)$ and
the nucleus of the loop it follows that $\alpha(a,x,y) = 1$ for
all $x, y \in Q \Leftrightarrow a \in N_{\lambda}(Q)$,
$\beta(a,x,y) = 1$ for all $x, y \in Q \Leftrightarrow a \in
N_{\rho}(Q)$ and $(a, N_{\lambda} \cap N_{\rho}) = 1 \Rightarrow a
\in N_{\mu}$. From here it follows that the equality (2.5) holds.

Particularly, if $Q$ is a Moufang loop, then from Lemma 2.4 it
follows that
$$a \in \mathcal{Z}(Q) \Leftrightarrow [a,x,y] = 1, [a,x] = 1
\quad \forall x, y \in Q. \eqno{(2.6)}$$ If $Q$ is an $IP$-loop or
$A$-loop, then from Lemmas 1,6 and 1.7 it follows that
$$a \in \mathcal{Z}(Q) \Leftrightarrow \alpha(a,x,y) = 1, (a,x) = 1
\quad \forall x, y \in Q \eqno{(2.7)}$$ or
$$a \in \mathcal{Z}(Q) \Leftrightarrow \beta(a,x,y) = 1, (a,x) = 1
\quad \forall x, y \in Q. \eqno{(2.8)}$$

\textbf{Lemma 2.7.} \textit{Let $H$ be a normal subloop of the
loop $Q$. Then the set $\mathcal{Z}_H(Q)$ is a normal subloop of
the loop $Q$ and $H \subseteq \mathcal{Z}_H(Q)$. Moreover, if $N$
is a subloop of $Q$ and $H \subseteq N \subseteq \mathcal{Z}_H(Q)$
then $N$ is a normal subloop of $Q$.} \vspace*{0.1cm}

\textbf{Proof.} Let $\varphi: Q \rightarrow \overline{Q} = Q/H$ be
the natural homomorphism. Denote by $\overline{A}$ the image of
set $A \subseteq Q$ under homomorphism $\varphi$. It is clear that
$A \subseteq H$ if and only if $\overline{A} = \{1\}$.
Particularly, the inclusion $\alpha(A,B,C) \subseteq H$ is
equivalent to the equality $\alpha(\overline{A}, \overline{B},
\overline{C}) = \{1\}$. Hence, due to one-to-one correspondence
between normal subloops of $Q$, containing $H$ and all normal
subloops of $\overline{Q}$ to prove the Lemma 2.7 it is sufficient
to consider that $H = E = \{1\}$. By (2.5) $\mathcal{Z}_E(Q)$
coincides with center $Z(Q)$, of loop $Q$, hence
$\mathcal{Z}_E(Q)$ is a normal subloop of $Q$. Then due to our
supposition  $\mathcal{Z}_H(Q)$ is also a normal subloop of $Q$.

Now let $N$ be a subloop of $Q$ and $N \subseteq
\mathcal{Z}_H(Q)$. Then $\alpha(N,Q,Q) \subseteq
\alpha(\mathcal{Z}_E(Q),Q,Q) = E \subseteq N$. Similarly,
$\beta(N,Q,Q) = N$, $(N,Q) = N$. By Lemma 2.5 $N$ is a normal
subloop of $Q$. $\Box$\vspace*{0.1cm}

Let $N$ be a normal subloop of the loop $Q$. We denote by
$\mathcal{A}^N(Q)$ the subloop of $Q$ generated by the set
$\alpha(N,Q,Q) \cup \beta(N,Q,Q) \cup(N,Q)$.  Particularly, if $Q$
is a  Moufang loop then $\mathcal{A}^N(Q)$ is generated by the set
\break $[N,Q,Q] \cup [N,Q]$.)\vspace*{0.1cm}

The subloop $\mathcal{A}^Q(Q)$ will be called
\textit{commutator-associator subloop} of the loop $Q$, and
sometimes is denoted by $Q^{(1)}$. For an arbitrary loop $Q$ the
subloop $\mathcal{A}^Q(Q)$ is generated by the set $\alpha(Q,Q,Q)
\cup  \beta(Q,Q,Q) \cup(Q,Q)$. But if $Q$ is a  Moufang loop then
$\mathcal{A}^Q(Q)$ is generated by the set $[Q,Q,Q] \cup [Q,Q]$
according to Lemma 2.4. \vspace*{0.1cm}

\textbf{Lemma 2.8.} \textit{Let $N$ be a normal subloop of the
loop $Q$. Then the subloop $\mathcal{A}^N(Q)$ will be a normal
subloop of the loop $Q$. Moreover, if $H$ is a subloop of $Q$ such
that $N \supseteq H \supseteq \mathcal{A}^N(Q)$ then $H$ will be a
normal subloop of the loop $Q$.} \vspace*{0.1cm}

\textbf{Proof.} The subloop $N$ is normal in $Q$, then by
Proposition 2.6 \break $\alpha(N,Q,Q) \subseteq N$, $\beta(N,Q,Q)
\subseteq N$, $(N,Q) \subseteq N$. Hence, $\mathcal{A}^N(Q)
\subseteq N$. From here it follows that for a subloop $H$ such
that $N \supseteq H \supseteq \mathcal{A}^N(Q)$ we have
$\alpha(H,Q,Q) \subseteq \alpha(N,Q,Q) \subseteq \mathcal{A}^N(Q)
\subseteq H$. Analogically, $\beta(H,Q,Q) \subseteq H$, $(H.Q)
\subseteq H$. Then by Proposition 2.6 the subloop $H$ is normal in
$Q$. $\Box$\vspace*{0.1cm}

\textbf{Corollary 2.9.} \textit{The commutator-associator subloop
$\mathcal{A}^Q(Q)$ of a loop $Q$ is the least normal subloop of
$Q$ such that the quotient loop $Q/\mathcal{A}^Q(Q)$ is an abelian
group.}\vspace*{0.1cm}

From the definitions of subloops $\mathcal{Z}_N$, $\mathcal{A}^N$
and Lemmas 2.7, 2.8 it follows.\vspace*{0.1cm}

\textbf{Corollary 2.10.} \textit{Let $N$ be a normal subloop of
the loop $Q$. Then the normal subloops $\mathcal{Z}_N(Q) =
\mathcal{Z}_N$}, $\mathcal{A}^N(Q) = \mathcal{A}^N$
\textit{satisfy the relations} $$\mathcal{Z}_N/N =
\mathcal{Z}(Q/N), \quad N/\mathcal{A}^N \subseteq
\mathcal{Z}(Q/\mathcal{A}^N),$$
$$\mathcal{Z}_{\mathcal{A}^N} \supseteq N, \quad
\mathcal{A}^{\mathcal{Z}_N} \subseteq N.$$

\textbf{Lemma 2.11.} \textit{Let $N$ be a normal subloop of the
loop $Q$ and let $H$ be the normal subloop defined in Lemma 2.5.
Then for any normal subloop $K$ of $Q$ $NK/K \subseteq Z(Q/K)$ if
and only if $H \subseteq K$.}\vspace*{0.1cm}

\textbf{Proof.} As $K$ is a normal subloop of the loop $Q$ then by
(1.4) $(xK)(yK \cdot nK) = (xK \cdot yK)(nK)$, $(nK \cdot xK)(yK)
= (nK)(xK \cdot yK)$, $nK \cdot xK = xK \cdot nK$  for all $n  \in
N$ and all $x, y \in Q$ when and only when $(x \cdot yn)K = xy
\cdot Kn$, $(nx \cdot y)K = nK \cdot xy$, $(nx)K =  x \cdot nK$
respectively. But this is equivalent to $x \cdot yn \in xy \cdot
Kn$, $nx \cdot y \in  nK \cdot xy$, $nx \in x \cdot nK$ or
$R^{-1}(n)L(x,y)n \in K$, $L^{-1}(n)R(x,y)n \in K$,
$L^{-1}(n)T(x)n \in K$. The centre $Z(Q/K)$ of the loop $Q/K$ is
an abelian group, then by (2.4) $NK/K \subseteq Z(Q/K)$ when and
only when $H \subseteq K$. $\Box$\vspace*{0.1cm}

Proceeding to the description of upper, lower central series. From
the definition of upper central series, (2.5) and Lemma 2.7 it
follows that: \vspace*{0.1cm}

\textbf{Proposition 2.12.} \textit{The transfinite upper central
series $\{\mathcal{Z}_{\alpha}\}$ of a loop $Q$ have the form:}

\textit{i) $\mathcal{Z}_0 = 1$;}

\textit{ii) for any ordinal $\alpha$  $\mathcal{Z}_{\alpha +
1}/\mathcal{Z}_{\alpha} =
\mathcal{Z}_{\alpha}(Q/\mathcal{Z}_{\alpha})$;}

\textit{iii) if $\alpha$ is a limit ordinal, $\mathcal{Z}_\alpha =
\bigcup_{\beta < \alpha} \mathcal{Z}_\beta$.}

\textit{Particularly, if $Q$ is a commutative loop either an
$IP$-loop, or a Moufang loop, or an $A$-loop then the normal
subloops $\mathcal{Z}_{\alpha}$ change in accordance with the
specifications shown after equalities(2.5) -- (2.8).}
\vspace*{0.1cm}

From the definitions of transfinite lower central series, Lemmas
2.8, 2.11 and Corollary 2.9 it immediately follows that:
\vspace*{0.1cm}

\textbf{Proposition 2.13.} \textit{The transfinite lower central
series $\{Q_{\xi}\}$ of a loop $Q$ have the form:}

\textit{$Q_0 = Q$, $Q_1 = \mathcal{A}^Q(Q)$, $Q_{\xi + 1} =
\mathcal{A}^{Q_{\xi}}(Q)$  for any ordinal $\xi$;}

\textit{$Q_{\xi} = \cap_{\eta < \xi}Q_{\eta}$ if $\xi$ is a limit
ordinal.}

\textit{If $Q$ is a Moufang loop  then the normal subloops change
in accordance with the definition of the subloop
$\mathcal{A}^N(Q)$, where $N$ is a normal subloop of
$Q$.}\vspace*{0.1cm}

Now we define the \textit{commutator-associator of type $(\alpha,
\beta)$ of weight $n$} inductively:

1) any associators of the form $\alpha(x,y,z)$, $\beta(x,y,z)$ and
any commutator $(x,y)$ are com\-mutator-associator of the type
$(\alpha, \beta)$ of weight $1$;

2) if $a$ is a commutator-associator of the type $(\alpha, \beta)$
of weight $n - 1$, then $\alpha(a,x,y)$, $\beta(a,x,y)$, $(a,x)$,
where $x, y \in Q$, are a commutator-associ\-ator of the type
$(\alpha, \beta)$ of the weight $n$.

If  only the associators of types  $\alpha$ and $\beta$
participate in the definition, then we get the \textit{associators
of type $(\alpha, \beta)$.}\vspace*{0.1cm}

We define by induction the \textit{commutator-associator of type
$(\mu)$ of weight $n$} :

1) any associator of form $[x,y,z]$ and any commutator $[x,y]$ are
commu\-ta\-tor-associators of type $(\mu)$ of weight $1$;

2) if $a$ is a commutator-associator of type $(\mu))$ of weight $n
- 1$, then $[a,x,y]$, $[a,x]$ are commutator-associators of type
$(\mu)$ of weight $n$.

If  only the associators of type $\mu$  participate in the
definition then we get the- \textit{associators of type
$(\mu)$}.\vspace*{0.1cm}

We assume that in the loop $Q$ all commutator-associators of
weight $n$ of one of the aforementioned types, for example
$(\alpha, \beta)$, are equal to unit $1$. Then we say that  the
loop $Q$ satisfies the \textit{commutator-associator identities}
of type $(\alpha, \beta)$ of weight $n$. The notions of
\textit{associator identities} of type $(\alpha, \beta)$ of weight
$n$, of \textit{associator identities} of type $(\alpha)$ of
weight $n$, \textit{associator identities} of type $(\beta)$ of
weight $n$, of \textit{associator identities} of type $(\mu)$ of
weight $n$ are  obvious.\vspace*{0.1cm}

Further we denote by $W_n(\alpha, \beta)$ the set of all
commutator-associators of type $(\alpha, \beta)$ of weight $n$, by
$W_n(\alpha)$ the set of all commutator-associators of type
$(\alpha)$ of weight $n$, by $W_n(\beta)$ the set of all
commutator-associators of type $(\beta)$ of weight $n$, by
$W_n(\mu)$ the set of all commutator-associators of type $(\mu)$
of weight $n$, and so on. Reinforce the Proposition 2.2 and
Proposition 2.4 in case of centrally nilpotent loops  using the
Propositions 2.12, 2.13.\vspace*{0.1cm}

\textbf{Theorem 2.14.} \textit{For a loop $Q$ the following
statements are equivalent:}

\textit{1) the loop $Q$ is centrally nilpotent of class $n$;}

\textit{2) the upper central nilpotent series of $Q$ have the
form}
$$E = \mathcal{Z}_0  \subset \mathcal{Z}_1 \subset \mathcal{Z}_2 \subset
\ldots \subset  \mathcal{Z}_{n-1} \subset \mathcal{Z}_n = Q,$$
\textit{where $E = \{1\}$, $\mathcal{Z}_{i+1}/\mathcal{Z}_i =
\mathcal{Z}(Q/\mathcal{Z}_i)$, $i = 0, \ldots, n - 1$,  and}
$$\mathcal{Z}_i = \{a \in Q \vert w(a, x_1, \ldots,
x_j) = 1  \forall w \in W_i(\alpha, \beta), \forall x_1, \ldots,
x_j \in Q\}. \eqno{(2.9)}$$

\textit{3) the lower central nilpotent series of $Q$ have the
form}
$$Q = \mathcal{A}_0 \supset \mathcal{A}_1 \supset \mathcal{A}_2
\supset \ldots \supset \mathcal{A}_{n-1} \supset \mathcal{A}_n =
\{1\},$$ \textit{where  $\mathcal{A}_{i + 1} =
\mathcal{A}^{\mathcal{A}_{i}}(Q)$ for $i = 0, 1, \ldots, n-1$ and
the normal subloop $\mathcal{A}_{i}$ coincides with the subloop
$A_i$ generated by all commutator-associators $w \in W(\alpha,
\beta)$.}

\textit{Particularly, it is sufficient to consider  $w \in
W_i(\mu)$ if $Q$ is a Moufang loop and $w \in W_i(\alpha)$ or $w
\in W_i(\beta)$ if $Q$ is an $IP$-loop or
$A$-loop.}\vspace*{0.1cm}

\textbf{Proof.} The equivalence of items 1) - 3) follows from
Proposition 2.2.

By definition $\mathcal{Z}_1 = \{a \in Q \vert \alpha(a,x,y) = 1,
\beta(a,x,y) = 1, (a,x) = 1 \forall x, y \in Q\}$. From
$\mathcal{Z}_{i+1}/\mathcal{Z}_i = \mathcal{Z}(Q/\mathcal{Z}_i)$
it follows that $\mathcal{Z}_{i+1}/\mathcal{Z}_i = \{a \in Q \vert
\alpha(a,x,y) \subseteq \mathcal{Z}_i, \beta(a,x,y) \subseteq
\mathcal{Z}_i, (a,x) \subseteq \mathcal{Z}_i \forall x, y \in
Q\}$. From here by induction it follows easily that (2.6).

From (2.9) it follows that all commutator-associators of type
$(\alpha, \beta)$ of weight $n - 1$ belong to the normal subloop
$\mathcal{Z}_1$, which according to (2.6) and Lemmas 1.6, 1.7
coincides with centre $Z(Q)$ of loop $Q$. Then the subloop
$A_{n-1}$, generated by these commutator-associators, is normal in
$Q$.

We assume that the subloop $A_{i+1}$ is normal in $Q$ and we
consider the natural homomorphism $\varphi: Q \rightarrow
Q/A_{i+1}$. The subloop $A_i/A_{i+1}$ belongs to the centre
$Z(Q/A_{i+1})$, then it is normal in $Q/A_{i+1}$. The inverse
image of $A_i/A_{i+1}$ under homomorphism $\varphi$ is $A_i$.
Hence, the subloop $A_i$ is normal in $Q$. $\Box$\vspace*{0.1cm}

\textbf{Corollary 2.15.} \textit{According to Proposition 2.13 let
$\{\mathcal{A}_{\xi}\}$, where $\mathcal{A}_{\xi + 1} =
\mathcal{A}^{\mathcal{A}_{\xi}}(Q)$, be the (transfinite) lower
central series of a loop $Q$. Then for any natural number $n$ the
normal subloop $\mathcal{A}_n$ of the series
$\{\mathcal{A}_{\xi}\}$ coincides with the subloop $A_n$ of the
loop $Q$, generated by all commutator-associators of type
$(\alpha, \beta)$ of weight $n$  and the quotient loop
$Q/\mathcal{A}_n$ is centrally nilpotent of class $\leq
n$.}\vspace*{0.1cm}

\textbf{Proof.} From Lemma 1.1, (1.4), (2.3) it follows that the
subloop $A_n$ is normal in $Q$.

By definition,  $\mathcal{A}_1 = \mathcal{A}^Q(Q)$ and from Lemma
2.8 it follows that $\mathcal{A}_1$ is the subloop of $Q$
generated by all commutator-associators of type $(\alpha, \beta)$
of weight 1. Hence $\mathcal{A}_1$ is the normal subloop of $Q$
generated by all commutator-associ\-ators of type $(\alpha,
\beta)$ of weight 1, i.e. $\mathcal{A}_1 = A_1$.

We consider the normal subloops $\mathcal{A}_{n+1}$ and $A_{n+1}$.
By construction, $\mathcal{A}_{n+1}$ is the subloop of $Q$
generated by the set $\alpha(\mathcal{A}_n,Q,Q) \bigcup
\beta(\mathcal{A}_n,Q,Q) \bigcup \break (\mathcal{A}_n,Q)$. By
inductive hypothesis $\mathcal{A}_n = A_n$.Then the set
$\alpha(\mathcal{A}_n,Q,Q) \bigcup \break \beta(\mathcal{A}_n,Q,Q)
\bigcup  (\mathcal{A}_n,Q)$ contain all commutator-associators of
type $(\alpha, \beta, 1)$ of weight $n + 1$. Hence
$\mathcal{A}_{n+1} \supseteq = A_{n+1}$ because
$\mathcal{A}_{n+1}$ is a normal subloop in $Q$.

Taking this into consideration, we consider the quotient loop
$Q/A_{n+1}$. In this loop all commutator-associators of type
$(\alpha, \beta)$ of weight $n+1$ are equal to unit. Then by (2.5)
all commutator-associators of weight $n$ will be in the centre of
the loop $Q/A_{n+1}$. Consequently, $\mathcal{A}_{n}/A_{n+1} =
A_n/A_{n+1} \subseteq \mathcal{Z}(Q/A_{n+1})$ as the centre of any
loop is a normal subloop. Proceeding to inverse images we get
$\alpha(\mathcal{A}_{n}, Q, Q), \beta(\mathcal{A}_{n}, Q, Q),
(\mathcal{A}_{n},Q) \subseteq A_{n+1}$. Hence \break
$\mathcal{A}_{n+1} \subseteq A_{n+1}$. Consequently, $A_{n+1} =
\mathcal{A}_{n+1}$.

From  $\mathcal{A}_{n} = A_n$ it follows that any
commutator-associators of type \break $(\alpha, \beta)$ of weight
$n$ of $Q/\mathcal{A}_{n}$ is  equal to unit. Then by Theorem 2.14
the loop $Q/\mathcal{A}_{n}$ is centrally nilpotent of class $\leq
n$. $\Box$\vspace*{0.1cm}

Reinforce the Corollary 2.15  in case of Moufang  loops and
$A$-loops.\vspace*{0.1cm}

\textbf{Proposition 2.16.} \textit{Let $Q$ be a Moufang loop
(respect. $IP$-loop or $A$-loop) and according to Proposition 2.13
let $\{\mathcal{A}_{\xi}\}$, where $\mathcal{A}_{\xi + 1} =
\mathcal{A}^{\mathcal{A}_{\xi}}(Q)$, be the (transfinite) lower
central series of a loop $Q$. Then for any natural number $n$ the
normal subloop $\mathcal{A}_n$ of the series
$\{\mathcal{A}_{\xi}\}$ coincides with the subloop $A_n$ of the
loop $Q$, generated by all commutator-associators of type $(\mu)$
(respect. of type $\alpha$ or $\beta$) of weight $n$  and the
quotient loop $Q/\mathcal{A}_n$ is centrally nilpotent of class
$\leq n$.}\vspace*{0.1cm}

\textbf{Proof.} It is easy to show by  induction on natural number
$n$ the equality $\mathcal A_n = A_n$ using the Lemmas 2.3 and
2.4. Then the Proposition 2.16 follows from Corollary
2.15.$\Box$\vspace*{0.1cm}

\textbf{Corollary 2.17.} \textit{A loop $Q$ is centrally nilpotent
of class $n$ when and only when it satisfies all
commutator-associators identities $w(x_1, \ldots, x_j) = 1$ of
type $(\alpha, \beta)$ of weight $n$,  $w \in W_n(\alpha, \beta)$,
but does not satisfy at least one identity $v(x_1, \ldots, x_l) =
1$ of type $(\alpha, \beta)$ of weight $n - 1$, $v \in
W_{n-1}(\alpha, \beta)$.}

\textit{Particularly, if $Q$ is a Moufang  loop, then it is
sufficient to consider $w \in W_n(\mu)$, $v \in W_{n-1}(\mu)$ and
if $Q$ is an $IP$-loop or $A$-loop, then it is sufficient to
consider $w \in W_n(\alpha$, $v \in W_{n-1}(\alpha)$ or $w \in
W_n(\beta)$, $v \in W_{n-1}(\beta)$ .}\vspace*{0.1cm}

The statements that follow from the equivalence of items 1), 3) of
Theorem 2.14, are true.\vspace*{0.1cm}

\section{On varieties  of Moufang loops  and  $A$-loops}

The notion of  loop $(Q, \cdot, 1)$ defined by one basic operation
(see (1.1))  is not algebraical. In order to apply the universal
algebraic techniques, one must use the universal algebraic
description of loops as algebras $(Q, \cdot, /, \backslash, 1)$
with  three binary operations, multiplication $\cdot$, left
division $\backslash$,  right division $/$ and one nullary
operation $1$, satisfying the identities
$$(x \cdot y)/y = x;
(x/y)\cdot y = x, \quad x\backslash(x\cdot y) = y, \quad
x\cdot(x\backslash y) = y, \eqno{(3.1)}$$
$$1x = x1 = x. \eqno{(3.2)}$$

\textbf{Proposition 3.1.} \textit{Every quasigroup (respect. loop)
$(Q, \cdot, /, \backslash)$ defined by  three basic binary
operations is a quasigroup (respect. loop) $(Q, \cdot)$ defined by
one basic binary operation}.\vspace*{0.1cm}

\textbf{Proof.} For quasigroup $(Q, \cdot, /, \backslash)$  we
consider the equation $ax = b$, where $a, b \in Q$. Then by (3.1)
$a \backslash (ax) = a \backslash b$, $x = a \backslash b$. If $ab
= ac$, where $a, b, c \in Q$, then $a \backslash (ab) = a
\backslash (ac)$, $b = c$. Hence the equation $ax = b$  has an
unique solution. Similarly it is proved that the equation $xa = b$
has an unique solution. Hence $(Q, \cdot)$ is a quasigroup with
one basic operation.  This completes the proof of Proposition
3.1.$\Box$\vspace*{0.1cm}

Let $(Q, \cdot, 1)$ be a Moufang loop. According to the definition
it  satisfies the equivalence identities (1.7) and the identities
(1.8). We denote $x \backslash y = x^{-1}y$,  $x / y = xy^{-1}$.
Then the algebra $(Q, \cdot, /, \backslash, 1)$ satisfies  the
identities (3.1) and (1.5). Hence $(Q, \cdot, /, \backslash, 1)$
is a Moufang loop. From here and Proposition 3.1 it follows that
for Moufang loops the definition of loop with three basic binary
operations and the definition of loop with one binary operation
coincide.\vspace*{0.1cm}

\textbf{Proposition 3.2.} \textit{The set $\frak M$ of all Moufang
loops defined by three basic binary operations form a variety
defined by one from equivalent identities (1.5).}\vspace*{0.1cm}

Obviously, from (3.2) it follows $x/x = x \backslash x = 1$, $1
\backslash x = x / 1 = x$ for every $x$ in $Q$, here $1$ is an
identity for $(Q, \cdot)$, a left identity for $Q, \backslash)$, a
right identity for $(Q, /)$. Assume that $1$ is a left identity
for $(Q, /)$, $1 / x = x$. Then by (3.1) $(1 / x)x = xx$, $1 = xx
= x^2$.$\Box$\vspace*{0.1cm}

We analyze the concept of $A$-loop. If for loop $(Q, \cdot, /,
\backslash, 1)$  we consider the group generalized by all left and
right translations with respect to basic operations $\cdot, /,
\backslash$ and  enter the concept of group of inner mappings,
then by above-stated we shall receive the identity $x^2 = 1$. The
$A$-loops with identity $x^2 = 1$ are investigated in \cite{GKN}.
We exclude this case. Moreover, in order not to limit the Bruck
and Paige's concept of $A$-loop $(Q, \cdot, 1)$ we
define.\vspace*{0.1cm}

Let $(Q, \cdot, /, \backslash, 1)$ be a loop, let $\frak{C}(Q)$ be
the multiplication group of loop $(Q, \cdot, 1)$ and let
$\frak{I}(Q)$ be the group of inner mapping of loop $(Q, \cdot,
1)$. The loop $(Q, \cdot, /, \backslash, 1)$ will called
\textit{$A$-loop} if the inner mappings in $\frak{I}(Q) \subseteq
\frak{C}(Q)$ are automorphism for loop $(Q, \cdot,
1)$.\vspace*{0.1cm}

Every  $A$-loop is power associative \cite{BP}, then $1 / x = x
\backslash 1$. Further for $A $-loops we shall use the designation
$$x^{-1} = 1 / x = x \backslash 1. \eqno{(3.3)}$$

For a loop $(Q, \cdot, /, \backslash, 1)$ we have $L(x)y = R(y)x =
xy$. Denote $L_{(/)}(x)y = R_{(/)}(y)x = x / y$,
$L_{(\backslash)}(x)y = R_{(\backslash)}(y)x = x \backslash y$.
From (3.1) it follows that $x / (y \backslash x) = y$, $(x / y)
\backslash x = y$, $L_{(/)}(x)R_{(\backslash)}(x)y = y$,
$R_{(\backslash)}(x)L_{(/)}(x)y = y$. Similarly, from (3.1) it
follows that $R_{(/)}R(y)(x) = x$, $L_{(\backslash)}L(x)y = y$.
Hence $$L^{-1}(x) = L_{(\backslash)}(x), \quad R^{-1}(x) =
R_{(/)}(x), \quad L_{(/)}(x) = R^{-1}_{(\backslash)}(x).
\eqno{(3.4)}$$

The group of inner mapping $\frak{I}(Q)$ of loop $(Q, \cdot, /,
\backslash, 1)$ is generated by inner mapping $T(x) =
L^{-1}(x)R(x), \quad R(x,y) = R^{-1}(xy)R(y)R(x), \quad L(x,y)
\break = L^{-1}(xy)L(x)L(y)$ by (1.1) or according to (3.4) is
generated by inner mapping $$\overline{T}(x) =
L_{(\backslash)}(x)R(x), \quad \overline{R}(x,y) =
R_{(/)}(xy)R(y)R(x),$$
$$\overline{L}(x,y) = L_{(\backslash)}(xy)L(x)L(y). \eqno{(3.5)}$$
From (3.4), (3.5) follows that $$\overline{T}^{-1}(x) =
R_{(/)}(x)L(x), \quad \overline{R}^{-1}(x,y) =
R_{(\backslash)}(x)R_{(\backslash)}(y)R(xy),$$
$$\overline{L}^{-1}(x,y) =
L_{(\backslash)}(y)L_{(\backslash)}(x)L(xy). \eqno{(3.6)}$$ The
group of inner mapping $\frak{I}(Q)$ of loop $(Q, \cdot, /,
\backslash, 1)$ is generated by inner mapping $T(x)$, $R(x,y)$, $
L(x,y)$. Then from (3.5), (3.6) it follows that any inner mapping
$\alpha \in \frak{I}(Q)$ has a representation $\alpha = S_1S_2
\ldots S_n$, where $$S_i \in \{\overline{T}(x_1),
\overline{R}(x_2,x_3), \overline{L}(x_4,x_5),
\overline{T}^{-1}(x_6), \overline{R}^{-1}(x_7,x_8),
\overline{L}^{-1}(x_9,x_{10})\},$$ $x_1, \ldots, x_{10} \in Q$. We
choose the least number $n$ with such a property and we fix such a
record. The number $n = n(\alpha)$ will called a \textit{length}
of mapping $\alpha$.

Let $(Q, \cdot, /, \backslash, 1)$ be an $A$-loop. Then
$\alpha(u\cdot u) = \alpha(u) \cdot \alpha(v)$ for all $u, v \in
Q$. We transform the last expression into loop expression
$\overline{\alpha}(u\cdot u) = \overline{\alpha}(u) \cdot
\overline{\alpha}(v)$ with respect to basic operations $\cdot, /,
\backslash$ with the help of relations $L_{(\backslash)}(x)y =
x\backslash y$, $R_{(/)}(x)y = y / x$. Further,
$\overline{\alpha}(u\cdot u) = \overline{\alpha}(u) \cdot
\overline{\alpha}(v)$ turns to identity if to consider $\cdot, /,
\backslash$ as symbols of operations and $x_1, \ldots, x_{10}, u,
v$ as free variables (see \cite[pag. 268]{Mal}). Obviously that
received identity of length $n$ is valid in  loop $(Q, \cdot, /,
\backslash, 1)$.

We denote by $\frak{A}_n$ the variety of all $A$-loops $(Q, \cdot,
/, \backslash, 1)$ defined by above described identities of length
$\leq n$. We get $\frak{A}_1 \supseteq \frak{A}_2 \supseteq \ldots
\frak{A}_n \supseteq \ldots$. Obviously, $\frak{A} = \cap
\frak{A}_n$ is a variety and is the variety of all $A$-loops with
three basic operations. Is easy to show induction on $n$ that
every identity $\overline{\alpha}(u\cdot u) = \overline{\alpha}(u)
\cdot \overline{\alpha}(v)$ of length $n$ is a consequence  of
identity of length $1$. Then $\frak{A} = \frak{A}_n$. Hence,
according to (3.5), (3.6), we proved.\vspace*{0.1cm}

\textbf{Proposition 3.3.}  \textit{The set $\frak A$ of all
$A$-loops defined by three basic binary operations form a variety
defined by six  identities} $S_i(u \cdot v) = S_iu \cdot S_iv$,
where $S_i \in \{\overline{T}(x_1), \overline{R}(x_2,x_3),
\overline{L}(x_4,x_5), \overline{T}^{-1}(x_6),
\overline{R}^{-1}(x_7,x_8),
\overline{L}^{-1}(x_9,x_{10})\}.$\vspace*{0.1cm}

 Now we follow Mal'cev \cite{Mal}.\vspace*{0.1cm}

1). Any variety of algebras $\frak B$ is quite determined by its
free algebra of countable infinite rank $F = f(\frak B)$ (Theorem
VI.13.3). Let $\frak C$ be a subvariety of $\frak B$. The least
congruence $\theta$ on $F$ such that $F/\theta \in \frak C$ is
called \textit{word congruence}.\vspace*{0.1cm}

2). The quotent algebra $F/\theta$ is a  free algebra of $\frak
C$.\vspace*{0.1cm}

3). A congruence $\theta$ on $F$ is word congruence if $x \theta
y$ implies $\varphi x \theta \varphi y$ for all endomorphisms
$\varphi; F \rightarrow F$. The congruence $\theta$ with such an
implication is called \textit{fully characteristic} (Corollary
VI.14.3).

4). Let $KI(F/\theta)$ denote the variety generated by $F/\theta$.
Then the mapping $\theta \rightarrow KI(F/\theta)$ is an
antiisomorphism of lattice of fully characteristic congruences and
lattice of all subvarieties of $\frak B$ (Corollary VI.14.4).

Further  we shall transfer the above-stated statements on a
variety $\frak L$ of all  loops defined by three basic binary
operations $\frak L$ and $\frak L$-free loop $F = F(\frak L$ with
infinite set of  free generators $f_1, f_2,
\ldots$.\vspace*{0.1cm}

5). According to \cite[pag. 61]{Mal} a relation of equivalence
$\theta$ of algebra $(A, \Omega)$ is a congruence if $\theta$ is
stable concerning all basic operations in $\Omega$.\vspace*{0.1cm}

6). If $(A, \Omega) \in \frak L$ then the conditions that $\theta$
is stable concerning all operations $(\cdot), (/), (\backslash)$:
$x_1 \theta x_2$, $y_1 \theta y_2$ $\rightarrow x_1 \alpha y_1
\theta x_2 \alpha y_2$, where $\alpha \in \{(\cdot), (/),
(\backslash)\}$, is equivalent to condition: $x \theta y
\Leftrightarrow xz \theta yz$, $zx \theta zy$ for all $z \in
Q$.\vspace*{0.1cm}

7). The relation $\theta$ with the least condition for loops $(Q,
\cdot, 1)$ in the literature \cite{Br58}, (\cite{Bel}) is called
\textit{normal congruence} and the inverse image $H$ of unity at
homomorphism $Q \rightarrow Q/\theta$ is called \textit{normal
subloop} of loop $(Q, \cdot, 1)$. It coincides with the notion of
normal subloop defined by (1.2).\vspace*{0.1cm}

8). From item 7) it follows that  all researches of previous
sections for loops defined by one basic binary operation literally
are valid for loops defined by three  basic binary
operations.\vspace*{0.1cm}

9).  According to items 3) and 7) we define. Let $\frak B$ be a
variety of loops defined by three basic binary operations and let
$L = L(\frak B)$ be a $\frak B$-free loop with infinite number of
free generators. A normal subloop $H$ of loop $L$ will be called
\textit{fully invariant} if $\varphi H \subseteq H$ for all
endomorphism $\varphi: L \rightarrow L$.\vspace*{0.1cm}

10). Let $X_{\infty}$ denote the free loop with infinite set of
free generators $x_1, x_2, \ldots$ of variety of all loops defined
by three basic binary operations. This free loop we shall use as
totality of "words": \textit{word} means an element in
$X_{\infty}$ in the alphabet $x_1, x_2, \ldots$.

If $Q$ is a loop and $\alpha$ is a mapping of free generators
$x_1, x_2, \ldots$ into $Q$ then the image of word $v = v(x_{i_1},
\ldots, x_{i_k}) \in X_{\infty}$ at homomorphism $\alpha:
X_{\infty} \rightarrow Q$ is called \textit{the value} of word $v$
in $Q$, $v(\alpha x_{i_1}, \ldots, \alpha x_{i_k})$.

The  word $v(x_{i_1}, \ldots, x_{i_k})$ is called an identity for
$Q$ if the unity $1$ is the  unique value of $v$  in $Q$, $v(a) =
1$ for all $a \in Q^k$.\vspace*{0.1cm}

11). Let $(Q, \cdot, /, \backslash, 1)$ be a loop, let $\varphi$
be an endomorphism of $(Q, \cdot, /, \backslash, 1)$, let
$\frak{I}(Q)$ be the group of inner mappings of loop $(Q, \cdot,
1)$ and let $\alpha \in \frak{I}(Q)$. Then $\varphi \frak{I}(Q) x
\subseteq \frak{I}(Q) \varphi x$ for all $x \in Q$.

Indeed, the inner mapping $\alpha \in \frak{I}(Q)$ is a product of
finite number of factors $T^{\pm 1}(a)$, $L^{\pm 1}(a, b)$,
$R^{\pm 1}(a, b)$, $a, b \in Q$. Then the assertion follows from
the  definition of endomorphism $\varphi$  of $(Q, \cdot, /,
\backslash, 1)$ and (3.4) -- (3.6).\vspace*{0.1cm}

12). Let $\frak U$ denote the variety of all loops defined by
three basic binary operations with $\frak U$-free loop $F =
f(\frak U)$ of infinite number of free generators $f_1, f_2,
\ldots$ and let  $\frak B$ be a subvariety of  $\frak U$ with
$\frak B$-free loop $G = G(\frak U)$ of infinite number of free
generators $g_1, g_2, \ldots$. A normal subloop $H$ of loop $Q \in
\frak B$ will be called \textit{a word subloop} if $H$ is  the
normal subloop of $Q$ generated as normal subloop by all values of
words of some given set  $M$: $M(Q)$ is the loop generated as
normal subloop of $Q$ by set $\{\varphi v \vert v \in M, \varphi
\in \hom(X_{\infty}, Q)\}$.\vspace*{0.1cm}

13). All words from $M$ are identities for a loop $Q \in \frak B$
if and only if  $M(Q) = \{1\}$. If $v \in M$ and $v = v(x_1,
\ldots, x_k)$, then $M(Q) = \{v(a_1, \ldots, a_k \vert \break
\forall a_1, \ldots, a_k \in Q\}$.\vspace*{0.1cm}

14). Let $M$ be some set of words. The word $u$ is said to be a
consequence of $M$ in a variety $\frak B$ if with  accuracy up to
renameng of variables   the word subloop generated by word $u$ is
a subloop of word subloop generated by set of word $M$ in a $\frak
B$-free loop. Two sets of words $M_1$ and $M_2$ are equivalent in
variety $\frak B$ if every word of $M_1$ is a consequence of $M_2$
and conversely.\vspace*{0.1cm}

15). A normal subloop $H$ of the loop $Q \in \frak B$ defined in
item 12) is a word subloop of $Q$ if and only if $H$ is fully
characteristic in $Q$.\vspace*{0.1cm}

The proof of this statement is the items 16) --
18).\vspace*{0.1cm}

16). Every word subloop $H$ of loop $Q \in \frak B$ is  fully
invariant in $Q$.\vspace*{0.1cm}

The assertion follows from the definition of word subloop, Lemma
1.3 and item 11).\vspace*{0.1cm}

17). Every fully invariant subloop $V$ of free loop $F = F(\frak
U)$ defined in item 12) is word subloop in $F$.\vspace*{0.1cm}

Indeed, the subloop $V$ is fully invariant in $F$. In such a case
it is possible to consider as  set  of words $M$ the set of such
words $v = v(x_{i_1}, \ldots, x_{i_k}) \in X_{\infty}$ that $v =
v(f_{i_1}, \ldots, f_{i_k}) \in V$.\vspace*{0.1cm}

18). Every fully invariant subloop of free loop $G = G(\frak B)$
defined in item 12) is word subloop in $G$.\vspace*{0.1cm}

Indeed, as $\frak B \subseteq \frak U$ then any   mapping
$\varphi$ of free generators $f_1, f_2, \ldots$ on free generators
$g_1, g_2, \ldots$ proceeds up to epimorphism $\varphi: F
\rightarrow G$ which results in representation $G \cong F/R$ with
property: $r(f) \in R$ for some $f = (f_{i_1}, \ldots, f_{i_k})$
implies  $r(\theta f) \in R$ at any endomorphism $\theta$ of loop
$F$. Hence $G$ can be  presented as quotient loop $G \cong F/R$ on
some of its fully characteristic subloop.

Let $N$ be the inverse image in $F$ of  given fully invariant
subloop of $G$. Then the subloop $N/R$ is fully invariant in
$F/R$. Every endomorphism of loop $F$ induces an endomorphism of
loop $G$, hence $\theta R \subseteq R$, $\theta N/R \subseteq
N/R$. From the definition of induced endomorphism it follows that
$\theta N \subseteq N$. Hence, the subloop $N$ is fully invariant
in $F$ and by item 17) a verbal subloop in $F$: $N = M(F)$ for
some set $M$ of works. But then $N/R = M(F/R)$, as
required.$\Box$\vspace*{0.1cm}

According to Propositions 3.2, 3.3 and item 8) we consider  the
variety $\frak U$ of all loops  either the variety $\frak M$ of
all Moufang loops, or the variety $\frak A$ of all $A$-loops
defined by three basic binary operations. Let $F$ denote the
$\frak U$-free either $\frak M$-free, or $\frak A$-free loop on a
countable infinite set of generators and let $\{F_{\xi}\}$ be the
(transfinite) lower central  series of  $F$. From Corollaries 2.15
and 2.17, Proposition 32.16 and above-stated items 1) -- 18) it
follows.\vspace*{0.1cm}

\textbf{Theorem 3.4}. \textit{Let $n$ be a finite  natural number
and let $F_n$ be  the  member of lower central series
$\{F_{\xi}\}$ of free loop $F$. Then}

\textit{1) $F_n$ is a fully invariant subloop of loop $F$
generated by all commutator-associators of type $(\alpha, \beta)$
of weight $n$.}

\textit{2) If $\frak N$ denotes the variety of all centrally
nilpotent loops  of class $\leq n$ in $\frak U$, then the quotient
loop $F/F_n$ is a free loop of variety $\frak N \subseteq \frak
U$.}

\textit{3) Let $W_n(\alpha, \beta)$ (respect. $W_n(\mu)$) denote
the set of all commutator-associ\-ators of type $(\alpha, \beta)$
of weight $n$ (respect.  the set of all commutator-associators of
type $(\mu)$ of weight $n$ either the set of all
commutator-associ\-ators of type $(\alpha)$ of weight $n$ or the
set of all commutator-associators of type $(\beta)$ of weight $n$)
of loop $X_{\infty}$. Then the variety $\frak N \supseteq \frak U$
is determined by the finite set of identities $w = 1$, where $w
\in W_n(\alpha, \beta)$. Particularly, if $\frak U$ is the variety
of all Moufang loops, then $w \in W_n(\mu)$ and if $\frak U$ is
the variety of all $A$-loops or all $IP$-loops, then $w \in
W_n(\alpha)$ or $w \in W_n(\beta)$.}

\section{Lower central series of Moufang loops  and \break $A$-loops}

According to  Proposition 2.13 and item 8) of Section 3 let
$$Q\supseteq \mathcal{A}_{1}(Q) \supseteq \ldots \supseteq
\mathcal{A}_{k}(Q) \supseteq \mathcal{A}_{k+1}(Q) \supseteq
\mathcal{A}_{k+2}(Q) \supseteq \ldots, \eqno{(4.1)}$$ where
$\mathcal{A}_{\xi + 1}(Q) = \mathcal{A}^{\mathcal{A}_{\xi}}(Q)$,
be the (transfinite) lower central series of a loop $(Q, \cdot, /,
\backslash, 1)$. Let $k$ be a finite natural number and let
$\mathcal{A}_i(Q) = \mathcal{A}_i$. By Corollary 2.10 and (2.5)
$$\mathcal{A}_{k+1}/\mathcal{A}_{k+2} \subseteq
\mathcal{Z}(Q/\mathcal{A}_{k+2}) = Z(Q/\mathcal{A}_{k+2}).
\eqno{(4.2)}$$ We consider the quotient loop $Q/\mathcal{A}_{k+2}
= G$ and we suppose that the image of series (4.1) under
homomorphism $\varphi: Q \rightarrow Q/\mathcal{A}_{k+2}$ has the
form
$$Q/\mathcal{A}_{k+2} \supset \mathcal{A}_{1}/\mathcal{A}_{k+2}
\supset \mathcal{A}_{k+1}/\mathcal{A}_{k+2} \supset
\mathcal{A}_{k+2}/\mathcal{A}_{k+2} = 1.$$ Let $a \in
\mathcal{A}_{k+2}$, $x, y \in Q$. Then by Corollary 2.15
$\alpha(a,x,y), \beta(a,x,y), (a,x) \in \mathcal{A}_{k+1}$ and by
(1.4) $\alpha(a,x,y)\mathcal{A}_{k+2},
\beta(a,x,y)\mathcal{A}_{k+2}, (a,x)\mathcal{A}_{k+2} \in
\mathcal{A}_{k+1}/\mathcal{A}_{k+2}$. From relation (4.2)
$\mathcal{A}_{k+1}/\mathcal{A}_{k+2} \subseteq
\mathcal{Z}(Q/\mathcal{A}_{k+2})$  we get
$$\alpha(a,x,y)\mathcal{A}_{k+2},
\beta(a,x,y)\mathcal{A}_{k+2}, (a,x)\mathcal{A}_{k+2} \in
\mathcal{Z}(Q/\mathcal{A}_{k+2}). \eqno{(4.3)}$$

According to Corollary 2.15 let $G \supset \mathcal{A}_1(G)
\supset \mathcal{A}_2(G) = e$ be the lower central series of loop
$G$.  If $H$ is a normal subloop of a loop $L$ then by (1.4)
$\alpha(xH, yH, zH) = \alpha(x,y,z)H$, $\beta(xH, yH, zH) =
\beta(x,y,z)H$, $(xH, yH) = (x,y)H$. By Corollary 2.15
$\mathcal{A}_t(L) = <w_t(x_1, \ldots, x_j) \vert \quad \forall w_t
\in W_t(\alpha, \beta) \break \forall x_1, \ldots, x_j \in L>$,
$\mathcal{Z}_t(L) = \{a \in L \vert w_t(a, x_1, \ldots, x_r) = 1
\quad \forall w_t \in W_t(\alpha, \beta) \break \forall x_1,
\ldots, x_r \in L\}$. From here it follows that
$$\mathcal{A}_{k+2}(Q) = \varphi^{-1}(e), \mathcal{A}_{k+1}(Q) \subseteq
\varphi^{-1}(\mathcal{A}_{1}(G)),  \mathcal{A}_{k}(Q) \subseteq G.
\eqno{(4.4)}$$

From (4.2),  (4.3)  it follows.\vspace*{0.1cm}

\textbf{Lemma 4.1.} \textit{Let $Q$ be a loop with (transfinite)
lower central series $\{\mathcal{A}_{\xi}(Q)\}$  and let $k$ be a
natural number. Then} $a \in \mathcal{A}_k(Q) \Rightarrow \break
\alpha(a,x,y)\mathcal{A}_{k+2}(Q),$
$\beta(a,x,y)\mathcal{A}_{k+2}(Q), (a,x)\mathcal{A}_{k+2}(Q) \in
Z(Q/\mathcal{A}_{k+2})(Q), \break  \forall x, y \in
Q$.\vspace*{0.1cm}

\textbf{Lemma 4.2.} \textit{Let $G$ be a centrally nilpotent loop
of class 2. Then \break $\alpha(a,x,y)$, $\beta(a,x,y)$, $(a,x)
\in Z(G)$ $\forall a, x, y \in G$, where $Z(G)$ means the centre
of loop $G$.}\vspace*{0.1cm}

The Lemma 4.2 follows from Lemma 4.1. It also easy follows from
Corollary 2.17 and (2.9).\vspace*{0.1cm}

From (4.2) it follows that the quotent loop
$\mathcal{A}_{k+1}/\mathcal{A}_{k+2}$ is an abelian group. Then
for $u, v \in \mathcal{A}_{k+1}$ we will write the expressions
\break $(u \backslash v) \pmod{\mathcal{A}_{k+2}}$,  $(u / v)
\pmod{\mathcal{A}_{k+2}}$ in the form $(u^{-1}v)
\pmod{\mathcal{A}_{k+2}}$, \break $(uv^{-1})
\pmod{\mathcal{A}_{k+2}}$ respectively.

The crucial result used in the proof of Lemma 4.3 is the Lemma
4.2. We will use it without mentioning it and for Moufang loops we
will use Lemma 4.2.\vspace*{0.1cm}

\textbf{Lemma 4.3.} \textit{Let  $G$ be a centrally nilpotent loop
of class 2 and let \break $a, b, x, y, z \in G$. If $G$ is a
Moufang loop then}
$$[a,x,y] = [x,y,a] = [y,a,x], \eqno{(4.5)}$$ $$[a,x,y]^{-1} = [a^{-1},x,y] = [a,x^{-1},y] =
[a,y,x], \eqno{(4.6)}$$
$$[ab, x] = [a,x][b,x][a,b,x]^3, [a,xy] = [a,x][a,y][a,x,y]^3, \eqno{(4.7)}$$
$$[a,xy,z] = [a,x,z][a,y,z], \quad [a, x,yz] =
[a,x,y][b,x,z], \eqno{(4.8)}$$
$$[ab, x,y] = [a,x,y][b,x,y]. \eqno{(4.9)}$$

\textit{If $G$ is an $A$-loop then} $$(ab, x) = (a,x)(b,x), (a,xy)
= (a,x)(a,y), \eqno{(4.10)}$$
$$\gamma(ab, x,y) = \gamma(a,x,y)\gamma(b,x,y),$$ $$\gamma(a,xy,z) =
\gamma(a,x,z)\gamma(a,y,z), \gamma(a, x,yz) =
\gamma(a,x,y)(b,x,z), \eqno{(4.11)}$$
$$(a\backslash b, x) = (a,x)^{-1}(b,x), (a,x\backslash y) =
(a,x)^{-1}(a,y), \eqno{(4.12)}$$
$$\gamma(a\backslash b, x,y) = \gamma(a,x,y)^{-1}\gamma(b,x,y),
\gamma(a,x\backslash y,z) =$$ $$\gamma(a,x,z)^{-1}\gamma(a,y,z),
\gamma(a, x,y\backslash z) = \gamma(a,x,y)^{-1}(b,x,z).
\eqno{(4.13)}$$
$$(a/b, x) = (a,x)(b,x)^{-1}, (a,x/y) = (a,x)(a,y)^{-1},
\eqno{(4.14)}$$
$$\gamma(a/b, x,y) = \gamma(a,x,y)\gamma(b,x,y)^{-1}, \gamma(a,x/y,z)
=$$ $$\gamma(a,x,z)\gamma(a,y,z)^{-1}, \gamma(a, x,y/z) \equiv
\gamma(a,x,y)(b,x,z)^{-1}, \eqno{(4.15)}$$ where $\gamma = \alpha$
or $\gamma = \beta$.\vspace*{0.1cm}

\textbf{Proof.} Let $G$ be a Moufang loop. The equalities (4.5),
(4.6) follow from items (iii), (iv), (1.9) of Lemma 1.4 with $x =
a$ in item (i), as $[a,x] \in Z(G)$, $[[a,x],y] = 1$.

The prime equality from (4.7) follows from (1.10)  with $x = a$,
$y = b$ in item (i) of Lemma 1.4. Any Moufang loop is an
$IP$-loop, then by (1.8) $(xy)^{-1} = y^{-1}x^{-1}$, $[x,y]^{-1} =
[y^{-1},x^{-1}]$, $[x,y,z]^{-1} = [z^{-1},y^{-1},x^{-1}]$. The
centre of any loop is an abelian group. Then the second equality
from (4.7) follows from item (ii) of Lemma 1.4 and (1.11) if used
$(xy)^{-1} = y^{-1}x^{-1}$ and replaced  $z^{-1}$ by $a$, $y^{-1}$
by $x$, $x^{-1}$ by $y$.

By (1.11), (1.12) $L(z,a)(xy)\cdot [a^{-1},z^{-1}] =
L(z,a)x(L(z,a)y \cdot [a^{-1},z^{-1}])$,
$$(xy)[xy,a,z]^{-1}\cdot [a^{-1},z^{-1}] =
((x[x,a,z]^{-1})(y[y,a,z]^{-1})) \cdot [a^{-1},z^{-1}].$$  As
$[a^{-1},z^{-1}] \in Z(G)$ then  by (1.11), (1.12)
$(xy)[xy,a,z]^{-1} = L(z,a)(xy) = L(z,a)x \cdot L(z,a)y =
(x[x,a,z]^{-1})(y[y,a,z]^{-1})$.  By (4.6) $[a,z,x], [a,z,y] \in
Z(G)$ implies $[x,a,z], [y,a,z] \in Z(G)$.  Then from
$(xy)[xy,a,z]^{-1}] = \break (x[x,a,z]^{-1})(y[y,a,z]^{-1})$ it
follows $[xy,a,z]^{-1} = [x,a,z]^{-1}[y,a,z]^{-1}$ which by (4.6)
implies (4.8).

The subloops $\mathcal{Z}_2$, $G_k$ are normal in $G$ and $a, b
\in \mathcal{Z}_2$ or $a, b \in G_k$. Then by (1.5) $L(y,x)a,
L(y,x)b \in \mathcal{Z}_2$ or $L(y,x)a, L(y,x)b \in
\mathcal{G}_k$. By (1.6), diassociativity of Moufang loops and
(4.6) $[x^{-1},y^{-1}] = x^{-1}y^{-1}xy$. We use (4.7), (4.8).
Then $[L(y,x)a, L(y,x)b, [x^{-1},y^{-1}]] = [L(y,x)a, L(y,x)b,
x][L(y,x)a, \break L(y,x)b, y][L(y,x)a, L(y,x)b, x]^{-1}[L(y,x)a,
L(y,x)b, y]^{-1} = 1$. According to (1.6) $L(y,x)a(L(y,x)b \cdot
[x^{-1},y^{-1}]) = (L(y,x)a \cdot (L(y,x)b)[x^{-1},y^{-1}])$.

By (1.12) $L(y,x)(ab)[x^{-1},y^{-1}] =  L(y,x)a(L(y,x)b \cdot
[x^{-1},y^{-1}])$. Then $L(y,x)(ab) = L(y,x)aL(y,x)b$. We use
(1.11). Then $(ab)[ab,x,y] = \break a[a,x,y] \cdot b[b,x,y]$.
Consequently, $[ab,x,y] = [a,x,y][b,x,y]$ and (4.9) is proved.

Now, let $G$ be a centrally nilpotent of class $2$  $A$-loop and
let $a, b, x, y, z \in G$. According to (2.6) and Corollary 2.17
we assume that $a, b \in Z_{2}(G)$. Then from (2.5) and
Corollary 2.17 it follows that the commutator-associators
$\alpha(a, x, y), \beta(a, x, y), (a, x)$ belong to center $Z(G)$
of loop $G$. By definition $ax \cdot y = a\alpha(a, x, y)\cdot xy,
y\cdot xa = yx\cdot \beta(a, x, y)a, ax = x(a(a, x))$. We prove
the identities (4.10) - (4.15) only for associators of type
$\alpha(a, x, y)$ as for $\beta(a, x, y), (a, x)$ the
corresponding identities are proved analogically.

As $\alpha(a, x, y)\in Z(G)$ then from $ax\cdot y = a\alpha(a, x,
y)\cdot xy$ we get $ax\cdot y = a(\alpha(a, x, y)x\cdot y)$,
$R(y)L(a)x = L(a)R(y)(\alpha(a, x, y)x)$, $S(a, y)x = \alpha(a, x,
y)x$, where $S(a, y) = R(y)^{-1}L(a)^{-1}R(y)L(a)$. Obviously,
$S(a, y)1 = 1$, i.e. $S(a, y)$ is an inner mapping. Then $S(a, y)$
is an automorphism of loop $G$. Hence $\alpha(a, xz, y)(xz) = S(a,
y)(xz) = S(a, y)x\cdot S(a, y)z = \alpha(a, x, y)x\cdot \alpha(a,
z, y)z = (\alpha(a, x, y)\alpha(a, z, y))(xz)$, i.e. $\alpha(a,
xz, y) = \alpha(a, x, y)\alpha(a, z, y))$. The identities
$\alpha(ab, x, y) = \alpha(a, x, y)\alpha(b, x, y)$, $\alpha(a, x,
yz) = \alpha(a, x, y)\alpha(a, x$, $z))$, $(ab, x) = (a, x)(b,
x)$, $(a, xy) = (a, x)(a, y)$ are proved by analogy. Consequently,
the identities (4.10), (4.11) hold for the associators of type
$\alpha(x, y, z)$ and the commutators $(x, y)$.

Further, according to (4.10) $(a, x) = (a/b \cdot b, x) = (a/b,
x)(b, x), (a/b, x) =$ $(a, x)(b, x)^{-1}$. The other identities
(4.10) - (4.15) are proved in a similar manner.
$\Box$\vspace*{0.1cm}

\textbf{Lemma 4.4.} \textit{Let  $Q$ be a loop with (transfinite)
lower central series $\{\mathcal{A}_{\xi}\}$, let $k$ be a natural
number and let $a, b \in \mathcal{A}_k(Q)$,  $x, y, z \in Q$.}

\textit{If $Q$ is a Moufang loop then}
$$[a,x,y] \equiv [x,y,a] \equiv [y,a,x] \pmod{\mathcal{A}_{k+2}},$$
$$[a,x,y]^{-1} \equiv  [a^{-1},x,y] \equiv  [a,x^{-1},y] \equiv
[a,y,x] \pmod{\mathcal{A}_{k+2}}, \eqno{(4.16)}$$
$$[ab, x] \equiv [a,x][b,x][a,b,x]^3 \pmod{\mathcal{A}_{k+2}},$$
$$[a,xy] \equiv [a,x][a,y][a,x,y]^3 \pmod{\mathcal{A}_{k+2}}, \eqno{(4.17)}$$
$$[a,xy,z] \equiv [a,x,z][a,y,z]\pmod{\mathcal{A}_{k+2}},$$ $$[a, x,yz] \equiv
[a,x,y][b,x,z] \pmod{\mathcal{A}_{k+2}}, \eqno{(4.18)}$$
$$[ab, x,y] \equiv [a,x,y][b,x,y] \pmod{\mathcal{A}_{k+2}}. \eqno{(4.19)}$$

\textit{If $G$ is an $A$-loop then} $$(ab, x) \equiv
(a,x)(b,x)\pmod{\mathcal{A}_{k+2}},$$ $$(a,xy) \equiv (a,x)(a,y)
\pmod{\mathcal{A}_{k+2}}, \eqno{(4.20)}$$
$$\gamma(ab, x,y) \equiv \gamma(a,x,y)\gamma(b,x,y)\pmod{\mathcal{A}_{k+2}},$$
$$\gamma(a,xy,z) \equiv
\gamma(a,x,z)\gamma(a,y,z) \pmod{\mathcal{A}_{k+2}},$$ $$
\gamma(a, x,yz) \equiv \gamma(a,x,y)(b,x,z)
\pmod{\mathcal{A}_{k+2}}, \eqno{(4.21)}$$
$$(a\backslash b, x) \equiv (a,x)^{-1}(b,x) \pmod{\mathcal{A}_{k+2}},$$
$$(a,x\backslash y) \equiv
(a,x)^{-1}(a,y) \pmod{\mathcal{A}_{k+2}}, \eqno{(4.22)}$$
$$\gamma(a\backslash b, x,y) \equiv
\gamma(a,x,y)^{-1}\gamma(b,x,y),$$ $$\gamma(a,x\backslash y,z)
\equiv \gamma(a,x,z)^{-1}\gamma(a,y,z) \pmod{\mathcal{A}_{k+2}},$$
$$ \gamma(a, x,y\backslash z) \equiv \gamma(a,x,y)^{-1}(b,x,z)
\pmod{\mathcal{A}_{k+2}}. \eqno{(4.23)}$$
$$(a/b, x) \equiv (a,x)(b,x)^{-1} \pmod{\mathcal{A}_{k+2}},$$ $$(a,x/y) \equiv
(a,x)(a,y)^{-1} \pmod{\mathcal{A}_{k+2}}, \eqno{(4.24)}$$
$$\gamma(a/b, x,y) \equiv \gamma(a,x,y)\gamma(b,x,y)^{-1} \pmod{\mathcal{A}_{k+2}},$$
$$\gamma(a,x/y,z) \equiv \gamma(a,x,z)\gamma(a,y,z)^{-1} \pmod{\mathcal{A}_{k+2}},$$
$$\gamma(a, x,y/z) \equiv \gamma(a,x,y)(b,x,z)^{-1}
\pmod{\mathcal{A}_{k+2}}, \eqno{(4.25)}$$ where $\gamma = \alpha$
or $\gamma = \beta$.\vspace*{0.1cm}

Lemma 6.4 follows from Lemma 4.3 and (4.4).\vspace*{0.1cm}

\textbf{Theorem 4.5.} \textit{Let the Moufang loop (respect.
$A$-loop) $Q$ with lower central series
$\{\mathcal{A}_{\alpha}(Q)\}$  be generated by the set $M$. Then
for any integer $n \geq 0$ the quotient loop
$\mathcal{A}_n(Q)/\mathcal{A}_{n+1}(Q)$ is an abelian group and is
generated by those cosets of $Q$ modulo $\mathcal{A}_{n+1}(Q)$
 that contain commutator-associators of type
$(\mu)$ (respect. $(\alpha, \beta)$) of weight $n$ of the elements
of $M$.}\vspace*{0.1cm}

\textit{Proof.} We use the induction with respect to $n$. For $n =
0$ the statement is obvious. We assume that the quotient loop
$\mathcal{A}_{n-1}(Q)/\mathcal{A}_{n}(Q)$ is generated by cosets
of $Q$ modulo $\mathcal{A}_{n+1}(Q)$ containing
commutator-associators of type $(\mu)$ (respect. $(\alpha,
\beta)$) of weight $n - 1$ of elements of $M$. By Corollary 2.153
$\mathcal{A}_{n}(Q)$ is generated by the elements $[a,x,y]$,
$[a,x]$ (respect. $\alpha(a,x,y)$, $\beta(a,x,y)$, $(a,x)$), where
$a \in \mathcal{A}_{n-1}(Q)$, $x, y \in Q$. By induction
hypothesis, $a = (a^{\pm 1}_1 \ldots a^{\pm 1}_r)\cdot z$, where
$z \in \mathcal{A}_{n}$, and each $a_i$ is a commutator-associator
of type $(\mu)$ (respect. $(\alpha,\beta)$) of weight $n-1$ of
elements of $M$. The elements $x, y$ are terms with respect to
variables $b_i, c_i \in M$. After using several times the
identities (4.16) -- (4.19) (respect. (4.20) -- (4.25)), we get
that $[a,x,y]$, $[a,x]$ (respect. $\alpha(a,x,y)$, $\beta(a,x,y)$,
$(a,x)$) are products of terms of the forms $[a_i,b_i,c_i]$,
$[a_i,b_i]$, $[z,x,y]$, $[z,x]$ (respect. $\alpha(a_i,b_i,c_i)$,
$\beta(a_i,b_i,c_i)$, $(a_i,b_i)$, $\alpha(z,x,y)$,
$\beta(z,x,y)$, $(z,x)$, where $b_i, c_i \in M$. Since $[z,x,y],
[z,x], \break \alpha(z,x,y), \beta(z,x,y), (z,x) \in
\mathcal{A}_{n+1}$ the proof end.$\Box$\vspace*{0.1cm}

\textbf{Corollary 4.6}. \textit{Any finitely generated centrally
nilpotent Moufang loop or $A$-loop $Q$ satisfies the maximum
condition for its subloops.}

\textit{Proof.} Let $Q = \mathcal{A}_{0} \supset \mathcal{A}_{1}
\supset \ldots \mathcal{A}_{n} = \{1\}$ be the lower central
series of the loop $Q$ and let $H$ be a subloop of $Q$. We denote
$H_i = H \cap \mathcal{A}_i$. Using (1.3) it is easy to see that
the subloop $H_i$ is normal in $H$. By homomorphism theorems we
get $H_i/H_{i+1} = H_i/(H \cap \mathcal{A}_{i+1}) \cong
H_i\mathcal{A}_{i+1}/\mathcal{A}_{i+1} \subseteq
\mathcal{A}_{i}/\mathcal{A}_{i+1}$. By Theorem 4.5
$\mathcal{A}_{i}/\mathcal{A}_{i+1}$ is a finitely generated
abelian group. Then $H_i/H_{i+1}$ also is finitely generated. From
here it follows that the subloop $H$ is finitely generated. $\Box$

\section{The basis of identities}

Let $(F, \cdot, /, \backslash, 1)$ (respect. $(F, \cdot, ^{-1},
1)$) be the free $A$-loop (respect. free Moufang loop) on a
countable infinite set of free generators, $g_0, g_1, g_2,
\ldots$. Follow Evans \cite{Evans}. For $A$-loops we use the
designation  (3.3), $x^{-1} = 1 / x = x \backslash 1$. A
\textit{simple associator} in $A$-loop (respect.  Moufang loop)
$F$ is defined as follows. Each $g_i$ and $g^{-1}_i$ is a simple
associator. If $u, v, w$ are simple associators in $A$-loop
(respect. Moufang loop) $F$, then so is $\alpha(u, v, w), \beta(u,
v, w)$ (respect. $[u,v,w]$).\vspace*{0.1cm}

A simple associator in $F$ is said to \textit{involve} the
generator $g_i$ if (i) it is $g_i$ or $g_i^{-1}$, (ii) if it is
$\alpha(u, v, w), \beta(u, v, w)$ for $A$-loop and $[u,v,w]$ for
Moufang loop where at least one of $u, v, w$ involves
$g_i$.\vspace*{0.1cm}

Note that from the definitions of associators $\alpha(a,b,c)$,
$\beta(a,b,c)$, $[a,b,c]$ and (1.4) it follows that if a simple
associator involves $g_i$, then it lies in the normal subloop of
$(F, \cdot, 1)$ generated by $g_i$.\vspace*{0.1cm}

For $i = 0, 1, 2, \ldots$ we define the endomorphisms $\delta_i$,
 of  $F$ given by $\delta_i g_i = 1$, $\delta_i
g_j = g_j$ for $i \neq j$. The kernel of $\delta_i$ is the normal
subloop of $F$ generated by $g_i$.

According to item 8) in Section 3 let $\{F_{\xi}\}$ be the lower
central series of loop $(F, \cdot, 1)$ (or,  that same, of loop
$(F, \cdot, /, \backslash, 1))$. Follow Evans \cite[Lemma
2]{Evans}.\vspace*{0.1cm}

\textbf{Lemma 5.1.} \textit{If $w$ is an element in $F$ and $w$
maps onto $1$ under the endomorphisms $\delta_i$, $i = 0, 1, 2,
\ldots, 2n$, then $w \in F_n$ can be written as a product of
simple associators each of which involves $g_0, g_1, \ldots,
g_{2n}$ and $w \in F_n$.}\vspace*{0.1cm}

\textbf{Proof.}  We use the Propositions 3.2, 3.3 and item 8) of
Section 3.

We prove the assertion of Lemma 5.1 by induction on $n$. For $n =
0$ we prove this for one generator $g_0$. The kernel of $\delta_0$
is the normal subloop in $F$ generated by $g_0$. According to
Lemma 1.1 $\ker \delta_0$  is generated as a subloop by set
$\{\varphi g_0 \vert \quad \varphi \in \frak I(Q)\}$, where $\frak
I(Q)$ is the inner mapping group of $Q$. By (1.11) and (2.3)
$$L(z.y)x = x[x,y,z]^{-1}, \quad R(b,c)a = a\alpha(a,b,c),
 L(c,b)a = \beta(a,b,c)a. \eqno{(5.1)}$$ The property of
being a product of simple associators involving $g_0$ is clearly
preserved under multiplication. Furthermore, from (5.1) it follows
that for the proof of Lemma 5.1 at $i = 0$ is sufficient to prove
that   the application of any inner mappings $R(u,v)$, $L(u,v)$ to
a product of simple associators involving $g_0$ also preserves
this property. For $A$-loops it follows from definition of
$A$-loop.

Let $F$ be a Moufang loop. We define $x \equiv 1$ if $x \in \ker
\delta_0$. From $x, y \in \ker \delta_0$ it follows  $[x, u, v]
\equiv 1$, $[y, u, v] \equiv 1$ for all $u, v \in F$ by (5.1). It
is necessary to show that $[xy, u, v] \equiv 1$. But it follows
from Lemmas 1.2, 1.3.

Clearly, that $\ker \delta_0 \subseteq F_0 = F$. Hence the lemma
is proved for $n = 0$.

We suppose  by the induction hypothesis that lemma is correct for
$n - 1$. By Theorem  4.5 $F_{n - 1}$ modulo $F_{n}$ is generated
by commutator-associator of the type $(\alpha, \beta)$ (or
commutator-associator of the type $(\mu)$) of weight $n - 1$ of
the variables $g_1, g_2, \ldots, $. Then $$w \cong \prod_{i \in I}
v_i^{\pm 1}\mod F_{n}, \eqno{(5.2)}$$ where the set $I$ is finite
and for all $i \in I$ the element  $v_i$ is a
commutator-associator of weight $n - 1$. As every
commutator-associator  of weight $n - 1$ contains at most $2n - 2$
variables, each $v_i$, $i \in I$, does not contain some variable
of the set $\{g_1, g_2, \ldots , g_{2n}\}$.

Let the commutator-associator $v_j$ from (5.2) contains the
variables \break $g_{j_1}, g_{j_2}, \ldots, g_{j_r}$. If $i \in
\{j_1, j_2, \ldots , j_r\}$, then $\delta_i v_j = 1$, in the
contrary case $\delta_i v_j = v_j$. Further, from Corollary 2.17
it follows that the subloop $F_n$ is fully invariant,
particularly, $\delta_i F_n \subseteq F_n$. To apply the
endomorphism $\delta_0$ at (5.2). We obtain
$$1 \cong \delta_i w  \cong \prod_{i \in I_i}
v_i^{\pm 1}\mod F_{n+1},$$ where the subset $I_0$ of $ I = \{0, 1,
2, \ldots, 2n\}$ consists of those $i$ for which $v_i$ does not
contain $g_0$. Hence we can rewrite $w$ as $$w  \cong \prod_{i \in
I \backslash I_0} v_i^{\pm 1}\mod F_{n}.$$ If $I = I_0$ then $w
\in F_{n}$ and the proof is complete. If $I \neq I_0$, then we
repeat the procedure above for the representation of $w$ and the
endomorphism $\delta_1$. In at most $2n$ steps we will obtain that
$w \in F_{n}$.$\Box$\vspace*{0.1cm}

\textbf{Lemma 5.2.} \textit{Let $w$ be a word in $F$ in $g_0, g_1,
\ldots , g_t$ where $t \geq 2n$. Then $$w = (\ldots (uv_0)v_1
\ldots )v_{t-1})v_t,$$ where $u \in F_n$ and each $v_i$ is a
product of words of the form $\delta_{i_s} \ldots
\delta_{i_2}\delta_{i_1} w^{\pm 1}$, for some nonempty sequence
$i_1, i_2, \ldots , i_s$ from $0, 1, 2, \ldots ,
t$.}\vspace*{0.1cm}

\textbf{Proof.} We follow Evans \cite[Lemma 3]{Evans}. Define
$\gamma_i$, $i = 0, 1, 2, \ldots$ on words in $F$ by $\gamma_i w =
w \cdot (\delta_i w)^{-1}$. If $F$ is the $A$-loop then $\gamma_i
w = w \cdot (\delta_i w)^{-1} = w \cdot (1 / \delta_i w) = w \cdot
\delta_i(1 / w) = w \cdot \delta_i w^{-1}$. Similarly $\delta_i w
= w \cdot \delta_i w^{-1}$ for Moufang loop $F$.

Let $u = \gamma_t \ldots \gamma_1\gamma_0 w$. Then $\delta_i w =
1$, for $i = 0, 1, 2, \ldots t$. By Lemma 5.1, $u \in F_n$. Now
$\gamma_0 w = w \cdot  \delta_0 w^{-1}$, $\gamma_1\gamma_0 w = (w
\cdot \delta_0 w^{-1})(\delta_1\delta_0 w \cdot \delta_1 w^{-1})$,
and in general $\gamma_t \ldots \gamma_1\gamma_0 w = (\ldots (w
w_0)  w_1  \ldots )  w_t$, where each $w_i$ is an expression  of
words of the form $\delta_{i_s}\ldots \delta_{i_2}\delta_{i_1}
w^{\pm 1}$ for some nonempty sequence $i_1, i_2, \ldots, i_s$ from
$0, 1, 2, \ldots t$. Hence $w = (\ldots ((uv_0)v_1) \ldots )v_t$,
where $v_0 = w_t^{-1}$, $v_1 = w_{t-1}^{-1}$, \ldots, $v_t =
w^{-1}_0$.$\Box$\vspace*{0.1cm}

\textbf{Lemma 5.3.} \textit{Any identity $w(x_1, x_2, \ldots ) =
1$ in a Moufang loop or $A$-loop defined by three basic binary
operations  is an equivalent to a finite collection of identities
$w_i = 1$, $i = 1, 2, 3, \ldots $, where some of the $w_i$ belong
to $F_n$ and the remainder involve at most $2n$
variables.}\vspace*{0.1cm}

\textbf{Proof.} If $w = 1$ involves fewer than $2n$ variables,
there is nothing to prove. If this is not the case, we use Lemma
5.2, writing $w = 1$ as
$$((\cdot (u \cdot v_1) \cdot v_2) \cdot \ldots) \cdot w_s = 1.$$
where  $v_i$ are products of words in $x_1, x_2, \ldots$ of the
form $\delta_{i_s} \ldots \delta_{i_2} \delta_1 w^{\pm 1}$.  It
follows that $w = 1$ is equivalent to $u = 1$ where $u \in F_n$,
and identities $\delta_{i_s} \ldots \delta_{i_2} \delta_1 w = 1$,
since each of the $\delta_{i_s} \ldots \delta_{i_2} \delta_1 w =
1$ is a consequence of $w = 1$. Now the $\delta_{i_s} \ldots
\delta_{i_2} \delta_1 w = 1$ all involve at least one less
variable than $w = 1$. We may repeat the process with these new
identities until we obtain, in addition to the identities in
$F_n$, identities involving at most $2n$ variables, as
required.$\Box$\vspace*{0.1cm}

\textbf{Theorem 5.4.} \textit{The identities of a centrally
nilpotent Moufang loop (respect. $A$-loop) of class $n$ are
finitely based. The basis consists from one associator identity of
type $(\mu)$ (respect. of type $(\alpha)$ or $(\beta)$) of weight
$n$ in $2n +1$ variables and a finite collection of identities in
no more than $2n$ variables.}\vspace*{0.1cm}

\textbf{Proof.} We use the Theorem 3.4. Let $Q$ be a centrally
nilpotent Moufang loop or $A$-loop of class $n$ and let  $(F,
\cdot, /, \backslash, 1)$ be the free Moufang loop or free $A$-lop
on a countable infinite set of free generators, $g_1, g_2,
\ldots$. Let $H$ be the subloop of $F$ generated by $g_1, g_2,
\ldots, g_{2n+1}$ and $W$ the word subloop generated by the
identities of $Q$.

By Lemma 5.3, $W$ is generated, as a fully invariant subloop of
$F$, by words in $g_1, g_2, \ldots, g_{2n}$ and words in $F_n$.
Since $Q$ is centrally nilpotent of class $n$, $W \supseteq F_n$
and so $W$ is generated, as a fully invariant subloop, by set of
words in $g_1, \ldots, g_{2n}$ and the single associator word of
type $(\mu)$ (respect. $(\alpha)$ or $(\beta)$)  of weight $n$ in
$2n +1$ variables.  By Corollary 4.6, the subloop of $H$ generated
by these words can be generated by a finite set of words. This
finite set of words generates $W$ as a fully invariant subloop of
$F$. Hence, $F/W$ is finitely based and so are the identities of
$Q$.$\Box$

Nicolae I. Sandu,

Tiraspol State University of Moldova,

Chisin\u{a}u, R. Moldova

sandumn@yahoo.com

\end{document}